\documentclass[11pt, epsf]{amsart}
\usepackage{epsfig}
\usepackage{eucal}
\usepackage{mathrsfs}
\usepackage{amssymb}
\usepackage{latexsym}
\usepackage{amsfonts}
\usepackage[latin1]{inputenc}

\newenvironment{namelist}[1]{%
\begin{list}{}
 {
   
   \settowidth{\labelwidth}{#1}
   \setlength{\leftmargin}{1.1\labelwidth}
  }
 }{%
\end{list}}
\newcommand{\bi}{\begin{namelist}}
\newcommand{\ei}{\end{namelist}}

\newcommand{\bb}{\mathbb}

\newcommand{\om}{\Omega}

\newtheorem{Def}{Definition}

\newtheorem{Theo}{Theorem}

\newtheorem{Lem}{Lemma}

\title{On G-convergence of positive Self-adjoint operators}
\author[Hasan Almanasreh \lowercase{and} Mahmoud Shalalfeh]{Hasan Almanasreh$^1$ \lowercase{and} Mahmoud Shalalfeh$^2$}
\thanks{1- Mathematics Department, Hebron University,
        P.O. Box 40, Hebron, West Bank, Palestine}
\thanks{\hspace{0.4cm}Department of Math. Sciences, University of Gothenburg,
        SE-412 96 Gothenburg, Sweden}
\thanks{\hspace{0.4cm}E-mail: hasanm@hebron.edu}
\thanks{2- Mathematics Department, Hebron University,
        P.O. Box 40, Hebron, West Bank, Palestine}
\thanks{\hspace{0.4cm}E-mail: mahmouds@hebron.edu}

\begin{document}
\maketitle
\begin{abstract} We apply G-convergence theory to study the asymptotic of the eigenvalue problems of positive definite bounded self-adjoint $h$-dependent operators as $h\to\infty$. Two operators are considered; a second order elliptic operator and a general linear operator. Using the definition of G-convergence of elliptic operator, we review convergence results of the elliptic eigenvalue problem as $h\to\infty$. Also employing the general definition of G-convergence of positive definite self-adjoint operator together with $\Gamma$-convergence of the associated quadratic form, we characterize the G-limit as $h\to\infty$ of the general operator with some classes of perturbations. As a consequence, we also prove the convergence of the corresponding spectrum.
\end{abstract}
\section{Introduction}

Heterogeneous structures of materials appear often in physics, chemistry, mechanics, life sciences, and engineering. Very often one is also led to consider heterogeneous structures with a very fine
and complicated microstructure. Phenomena like heat conduction or transport equation are such structures which typically modeled by mathematical systems such as ordinary differential equations or partial differential equations, where the presence of fine microscopic scale is reflected in rapid oscillations of the coefficients. This situation can not in general be treated directly, and if it could be, the numerical methods employed to solve the problem require very fine degree of resolution in order the mesh can capture the oscillations, this of course on the count of costs, and in some situations, despite of mesh refinement, the solution will be out of reach.

In the present work we consider the convergence of the eigenvalue problem
\vspace{-0.2cm}
\begin{equation*}
\mathbf{H}_hu_h=\lambda_hu_h
\end{equation*}
for two different operators $\mathbf{H}_h$ in some suitable Hilbert spaces. We first consider $\mathbf{H}_h=-{\rm div}(\text{\textbf{A}}_h(x)\nabla)$ defined on $L^2(\Omega)$ with domain $H^1_0(\Omega)$, where $\textbf{A}_h(x)$ is some admissible coefficient matrix, and $\Omega$ is an open bounded subset of $\mathbb{R}^N$, $N\geq1$. Then we consider another arbitrary operator $\mathbf{H}_h=\mathbf{H}_0+V_h$ defined on $L^2(\Omega)$, $\Omega$ is as before, $\mathbf{H}_0$ is a positive definite bounded self-adjoint operator, and $V_h$ is a positive bounded Hermitian multiplicative perturbation. We are interested in the behavior of the operator $\mathbf{H}_h$ as the parameter $h\to\infty$, particularly we are interested in the asymptotic behavior of the point spectrum (the eigenvalues).

We will use classical operator and variational convergence theory. G-convergence theory is well-known for its applications in homogenization of partial differential equations. The concept was introduced in the late 1960's \cite{DEGS,SPA67,SPA68,SPA75} for linear elliptic and parabolic problems with symmetric coefficient matrices. Then the concept was extended to non-symmetric coefficient matrices \cite{MUR,TAR1,TAR2,TAR3} known as H-convergence. The definition was then generalized to positive definite self-adjoint operators \cite{DAL}. Later on, plenty of invaluable results are achieved for the elliptic and hyperbolic problems. In \cite{BRAC,CHDA,CHDE} G-convergence of monotone operators is proved. In \cite{SVA99,SVA00,SVA05} G-convergence of nonlinear parabolic operators is studied. The theory of G-convergence of differential operators in general is treated in \cite{ZKO,ZKOk}. Through out this paper, we will use the name G-convergence of the case of non-symmetric matrices as well.

The study of G-convergence of operators is often associated to the study of the asymptotic behavior for the associated quadratic forms in the calculus of variations via the notion of $\Gamma$-convergence which was introduced in the mid 1970's \cite{DEGF}. Here, we utilize and combine the two concepts in order to prove G-compactness for the operator $\mathbf{H}_h=\mathbf{H}_0+V_h$.

For the operator $\mathbf{H}_h=-{\rm div}(\text{\textbf{A}}_h(x)\nabla)$, the coefficient matrix $\text{\textbf{A}}_h$ is positive definite and bounded, then by the G-compactness criterion for elliptic operators, $\mathbf{H}_h$ has a G-limit as $h\to\infty$. The operator $\mathbf{H}_h=\mathbf{H}_0+V_h$ is positive definite, bounded, and self-adjoint, then using $\Gamma$-convergence for its associated quadratic form and the relation between G-convergence and $\Gamma$-convergence, we prove that $\mathbf{H}_h$ admits a G-limit as $h\to\infty$. Under suitable assumptions on the coefficient matrix $\textbf{A}_h(x)$ and on the perturbation $V_h$ we characterize the G-limits. Consequently we prove the convergence of the corresponding eigenvalues.\\

The paper is arranged as follows: In Section 2 we provide the reader with basic preliminaries on G-convergence and $\Gamma$-convergence. In Section 3 we revisit G-convergence of elliptic operators, and study the convergence properties of the corresponding eigenvalue problems. In Section 4 we prove the G-limit of the operator $\mathbf{H}_h=\mathbf{H}_0+V_h$.
\section{Preliminaries}
In what follows $\Omega$ will be an open bounded subset of $\mathbb{R}^N$, $N\geq1$, and the notations $\rightharpoonup$ and $\ast\negthickspace\rightharpoonup$ will denote weak and weak$^\ast$ convergence respectively. The domain is denoted by $\mathbf{D}$. Also $c$ and $C$ will denote real constants that might be different at each occurrence and are independent of all parameters, unless otherwise explicitly specified. The scalar products and norms are denoted by $\langle\cdot,\cdot\rangle$ and $\|\cdot\|$ respectively, where the norms $\|\cdot\|$ will be given indices to distinguish between them, while $\langle\cdot,\cdot\rangle$ are left without indices and their current definitions are obvious from the content.
\subsection{G-convergence}

For comprehensive materials on G-convergence we refer to e.g. \cite{DEF,OLS,SIL}, and for a general setting to positive definite self-adjoint operators to the monograph \cite{DAL}. Below we state two definitions of G-convergence; of elliptic operators and the general definition of positive definite self-adjoint operators.\\

Consider two positive real numbers $\alpha$ and $\beta$ such that $0<\alpha\leq\beta<\infty$, and define the following set of matrices\\
$
S(\alpha,\beta,\Omega)=\{\text{\textbf{A}}\in L^\infty(\Omega)^{N\times N}\,;\; (\text{\textbf{A}}(x,\xi),\xi)\geq\alpha |\xi|^2\,\text{and}\, |\text{\textbf{A}}(x,\xi)|\leq\beta|\xi|\; , \forall\xi\in\mathbb{R}^N$ and a.e $\;x\in\Omega\}\,.$\\
We shall define G-convergence of the following sequence of elliptic Dirichlet boundary value problem
\begin{equation}\label{16}
\left\{ \begin{array}{l}
-{\rm div}(\text{\textbf{A}}_h(x,Du_h)) = f \mbox{ in }\om,\\
u_h\in H_0^1(\Omega).
\end{array} \right.
\end{equation}

\begin{Def}\emph{ The sequence $\text{\textbf{A}}_h$ in $S(\alpha,\beta,\Omega)$
is said to be $G$-convergent to $\text{\textbf{A}}\in S(\alpha,\beta,\Omega)$, denoted as $\text{\textbf{A}}_h\xrightarrow{\;\,\text{\tiny G}\;\,}\text{\textbf{A}}$, if for every $f\in H^{-1}(\Omega)$, the sequence $u_h$ of solutions  of (\ref{16}) satisfies
$$
\left. \begin{array}{l}
u_{h} \rightharpoonup {u}\;\mbox{ in }\;H_0^1(\Omega), \\
\text{\textbf{A}}_{h}(\cdot,Du_{h})\;\rightharpoonup {\text{\textbf{A}}(\cdot,Du)}\;\mbox{ in }
\;[L^2(\Omega)]^N,
\end{array} \right.
$$
where $u$ is the unique solution of the problem
\begin{equation}\label{17}
\left\{ \begin{array}{l}
-{\rm div}(\text{\textbf{A}}(x,Du)) = f \mbox{ in }\Omega,\\
u\in H_0^1(\Omega).
\end{array} \right.
\end{equation}
}
\end{Def}
In the sequel we will only consider the case of linear coefficients matrix $\text{\textbf{A}}_h$, i.e., from now on  $\text{\textbf{A}}_h(x,\xi) =\text{\textbf{A}}_h(x)\xi$.

Here are some results that will be used later. These results are given without proofs, for the proofs we refer to \cite{DEF,MUR}.
\begin{Theo} G-compactness Theorem.

\emph{
For every sequence $\text{\textbf{A}}_h$ in $S(\alpha,\beta,\Omega)$ there exists a subsequence, still denoted by $\text{\textbf{A}}_h$, and a map $\text{\textbf{A}}\in S(\alpha,\beta,\Omega)$, such that
$\text{\textbf{A}}_h\xrightarrow{\;\,\text{\tiny G}\;\,}\text{\textbf{A}}$.
}
\end{Theo}
\begin{Theo} Uniqueness and Locality of G-limit.
\emph{
\begin{itemize}
\item [$(i)$] $\text{\textbf{A}}_h$ has at most one G-limit.
\item [$(ii)$] If $\text{\textbf{A}}_h=\tilde{A}_h$ on $\omega\subset\subset\Omega$ and $\text{\textbf{A}}_h\xrightarrow{\;\,\text{\tiny G}\;\,}\text{\textbf{A}}$ and $\tilde{\text{\textbf{A}}}_h \xrightarrow{\;\,\text{\tiny G}\;\,}\tilde{\text{\textbf{A}}}$ then $\text{\textbf{A}}=\tilde{\text{\textbf{A}}}$ on $\omega$.
\end{itemize}
}
\end{Theo}
\begin{Theo}
\emph{
If $\text{\textbf{A}}_h\xrightarrow{\;\,\text{\tiny G}\;\,}\text{\textbf{A}}$, then $\text{\textbf{A}}_h^t\xrightarrow{\;\,\text{\tiny G}\;\,}\text{\textbf{A}}^t$.
}
\end{Theo}
Let $\mathcal{Y}$ be a Hilbert space, we provide below the general definition of G-convergence, first we set some useful definitions.
\begin{Def}
\emph{
A function $F:\mathcal{Y}\to[0,\infty]$ is said to be lower semi-continuous ($lsc$) at $u\in\mathcal{Y}$, if $$F(u)\leq \sup_{U\in\mathcal{N}(u)}\inf_{v\in U}F(v)\, ,$$ where $\mathcal{N}(u)$ is the set of all open neighborhoods of $u$ in $\mathcal{Y}$.
}
\end{Def}

As a consequence of the above definition we have the following
\begin{itemize}
\item [$(i)$] The inequality in the above definition can be replaced by equality due to the fact that $F(u)\geq\inf \{F(v), v\in U\},\;\; \forall U\!\in\mathcal{N}(u)$.
\item [$(ii)$] $F$ is $lsc$ on $\mathcal{Y}$, if it is so at each $u\in\mathcal{Y}$.
\end{itemize}

\begin{Def}
\emph{
A function $F$ in $\mathcal{Y}$ is called quadratic form if there exists a linear dense subspace $\mathcal{X}$ of $\mathcal{Y}$
and a symmetric bilinear form $B:\mathcal{X}\times \mathcal{X}\to[0,\infty)$ such that
$$
F(u)=\left\{ \begin{array}{ll}
B(u,u)\, ,&\forall u\in \mathcal{X}\, ,\\
\infty\, ,&\forall u\in\mathcal{Y}\backslash \mathcal{X}.
\end{array}\right.
$$
}
\end{Def}

The following result provides a useful criterion for G-convergence of positive definite self-adjoint operators. See \cite{DAL} for the proof.
\begin{Lem}
\emph{
Given $\lambda>0$, $A_h\in\mathcal{P}_\lambda(\mathcal{Y})$, and an orthogonal projection $P_h$ onto $\mathscr{V}_h$. Suppose that for every $u\in\mathcal{Y}$, $A_h^{-1}P_hu$ converges in $\mathcal{Y}$, then there exists an operator $A\in\mathcal{P}_\lambda(\mathcal{Y})$ such that $A_h\xrightarrow{\;\,\text{\tiny G}\;\,}A$ in $\mathcal{Y}$.
}
\end{Lem}

\subsection{$\Gamma$-convergence}
For comprehensive introductions to $\Gamma$-convergence we refer to the monographs \cite{BRA, DAL}.\\
Let $\mathcal{Y}$ be a topological space, and let $F_h$ be a sequence
of functionals from $\mathcal{Y}$ to $\overline{\mathbb{R}}=\mathbb{R}\cup\{-\infty,\infty\}$.
\begin{Def}\emph{ A sequence of functionals $F_h:\mathcal{Y}\to \overline{\mathbb{R}}$ is said to be $\Gamma$-convergent to $F:\mathcal{Y}\to \overline{\mathbb{R}}$, written as $F(u)=\Gamma-\displaystyle\lim_{h\to\infty}F_h(u)$ and denoted by $F_h\xrightarrow{\;\,{\tiny \Gamma}\;\,}F$ if
 $$
 F(u)=\Gamma-\displaystyle\liminf_{h\to\infty}F_h(u)=\Gamma-\displaystyle\limsup_{h\to\infty} F_h(u)\,,
 $$
where $\Gamma-\displaystyle\liminf_{h\to\infty}$ and $\Gamma-\displaystyle\limsup_{h\to\infty}$ are the $\Gamma$-lower and $\Gamma$-upper limits respectively defined by
$$
F^i(u):=\Gamma-\liminf_{h\to\infty} F_h(u)=\sup_{U\in{\mathcal N}(u)}\liminf_{h\to\infty}
\inf_{v\in U}F_h(v)
$$
and
$$
F^s(u):=\Gamma-\limsup_{h\to\infty} F_h(u)=\sup_{U\in{\mathcal N}(u)}\limsup_{h\to\infty}
\inf_{v\in U}F_h(v).
$$
}
\end{Def}

$\Gamma$-convergence possesses the compactness property, that is, if $\mathcal{Y}$ is a separable metric space, then every sequence of functionals $F_h:\mathcal{Y}\to \overline{\mathbb{R}}$ has a $\Gamma$-convergent subsequence.\\

The following theorem is the cornerstone of the relation between $\Gamma$-convergence of quadratic forms of the class $\mathcal{Q}_\lambda(\mathcal{Y})$ (respectively $\mathcal{Q}_0(\mathcal{Y})$) and G-convergence of the associated operators of the class $\mathcal{P}_\lambda(\mathcal{Y})$ for $\lambda>0$ (respectively strong resolvent convergence of the associated operators of the class $\mathcal{P}_0(\mathcal{Y})$ ). For the proof of this theorem we refer to \cite{BRA, DAL}.
\begin{Theo}
\emph{
Let $\lambda>0$ be a real number, $F_h$ and $F$ be elements of $\mathcal{Q}_0(\mathcal{Y})$, and let $A_h\,,\,A\in\mathcal{P}_0(\mathcal{Y})$ be the associated operators respectively. Then the following are equivalent
\begin{itemize}
\item [(a)] $F_h\xrightarrow{\;\,{\tiny \Gamma}\;\,}F$.
\item [(b)] $(F_h+\lambda||\cdot||^2_\mathcal{Y})\xrightarrow[]{\;\,{\tiny \Gamma}\;\,}(F+\lambda||\cdot||^2_\mathcal{Y})$.
\item [(c)] $(A_h+\lambda I)\xrightarrow{\;\,{\tiny G}\;\,}(A+\lambda I)$.
\item [(d)] $A_h\to A$ in the SRS.
\end{itemize}
Also if $F_h\,,\,F\in\mathcal{Q}_\mu(\mathcal{Y})$ for $\mu>0$, and $A_h\,,\,A\in\mathcal{P}_\mu(\mathcal{Y})$ are the associated operators respectively, then the following are equivalent
\begin{itemize}
\item [(e)] $F_h\xrightarrow{\;\,{\tiny \Gamma}\;\,}F$.
\item [(f)] $A_h\xrightarrow{\;\,{\tiny G}\;\,}A$.
\end{itemize}
}
\end{Theo}

\section{G-convergence of elliptic operators}

In this section we review some basic results of G-convergence of elliptic operators with source function $f_h$, where the main task is the discussion of eigenvalue problems ($f_h=\lambda_hu_h$). Before proceeding, a time is devoted to study the Dirichlet boundary value problem with $h$-dependent source function, which turns out to be useful in setting the results of the corresponding eigenvalue problem. The following two lemmas are required in proving the below homogenization results, we refer to \cite{MUR} for the proofs.
\begin{Lem}
\emph{
Let $\xi_h\in [L^2(\Omega)]^N$ be weakly convergent to $\xi$ in $[L^2(\Omega)]^N$, and $u_h\in H^1(\Omega)$ weakly convergent to $u$ in $H^1(\Omega)$, if
$$
{\rm div}(\xi_h)\rightarrow {\rm div}(\xi)\;\;\text{in}\; H^{-1}(\Omega)\, ,
$$
then
$$
\langle\xi_h,Du_h\rangle\;\ast\negthickspace\rightharpoonup \langle\xi,Du\rangle\;\;\text{in}\;D^\star(\Omega)\, ,
$$
where $D^\star(\Omega)$ is the dual space of the dense space $D(\Omega)=C^\infty_0(\Omega)$.
}
\end{Lem}

Before proving Theorem 5, we state and prove the following theorem for elliptic boundary value problem with source function $f_h$.
\begin{Theo}
\emph{
For the Dirichlet boundary value problem
\begin{equation}\label{34}
\left\{ \begin{array}{ll}
-{\rm div}(\text{\textbf{A}}_h(x)\nabla u_h) =f_h & \text{in}\; \Omega\, , \\
u_h\in H_0^1(\Omega),
\end{array} \right.
\end{equation}
if $\text{\textbf{A}}_h\in S(\alpha,\beta,\Omega)$ and if $f_h$ converges in $H^{-1}(\Omega)$ to $f$, then the sequence $u_h$ of solutions to (\ref{34}) is weakly convergent in $H_0^1(\Omega)$ to the solution of
\begin{equation}\label{34_2}
\left\{ \begin{array}{ll}
-{\rm div}(\text{\textbf{A}}(x)\nabla u) =f & \text{in}\; \Omega\, , \\
u\in H_0^1(\Omega),
\end{array} \right.
\end{equation}
where $\text{\textbf{A}}$ is the G-limit of $\text{\textbf{A}}_h$.
}
\end{Theo}
\hspace{-4mm}\emph{\underline{Proof}}. The weak form of (\ref{34}) is to find $u_h\in H_0^1(\Omega)$ such that $\forall v\in H_0^1(\Omega)$
\begin{equation}\label{35}
a_h(u_h,v)=\langle f_h,v\rangle,
\end{equation}
where $a_h(u_h,v)=\langle \text{\textbf{A}}_h \nabla u_h,\nabla v\rangle$. Since $\text{\textbf{A}}_h\in S(\alpha,\beta,\Omega)$, we have the following a priori estimate
$$
\alpha||u_h||_{H_0^1(\Omega)}^2\leq a_h(u_h,u_h)=\langle f_h,u_h\rangle\leq c||f_h||_{H^{-1}(\Omega)}||u_h||_{H_0^1(\Omega)}\, ,
$$
hence
\begin{equation}\label{36}
||u_h||_{H_0^1(\Omega)}\leq \frac{C}{\alpha}.
\end{equation}
By (\ref{36}) and the upper bound of $\text{\textbf{A}}_h$
\begin{equation}\label{37}
||\text{\textbf{A}}_h\nabla u_h||_{L^2(\Omega)}\leq C\frac{\beta}{\alpha}.
\end{equation}
So both $u_h$ and $\text{\textbf{A}}_h\nabla u_h$ are bounded sequences in $H_0^1(\Omega)$ and $[L^2(\Omega)]^N$ respectively, therefore, up to subsequences still denoted by $u_h$ and  $\text{\textbf{A}}_h\nabla u_h$
\begin{equation}\label{38}
u_h\rightharpoonup u\;\;\text{in}\; H_0^1(\Omega)
\end{equation}
and
\begin{equation}\label{39}
\text{\textbf{A}}_h\nabla u_h\rightharpoonup \mathcal{M}\;\;\text{in}\; [L^2(\Omega)]^N\, .
\end{equation}
\textbf{Claim}: we argue that $\mathcal{M}=\text{\textbf{A}}\nabla u$, where $\text{\textbf{A}}$ is the G-limit of $\text{\textbf{A}}_h$ (the existence and uniqueness of $\text{\textbf{A}}$ is guaranteed by virtue of Theorems 1 and 2).\\
\underline{\emph{Proof of the claim}}. By (\ref{39}) it holds that
\begin{equation}\label{40}
-{\rm div}(\text{\textbf{A}}_h\nabla u_h)\rightharpoonup -{\rm div}(\mathcal{M})\;\;\text{in}\; H^{-1}(\Omega)\,,
\end{equation}
which means that $\forall v\in H_0^1(\Omega)$
\begin{equation}\label{41}
\lim_{h\to\infty}\langle-{\rm div}(\text{\textbf{A}}_h\nabla u_h)\, ,\, v\rangle = \langle-{\rm div}(\mathcal{M})\, ,\, v\rangle.
\end{equation}
Since $f_h$ converges to $f$ in $H^{-1}(\Omega)$ and by (\ref{34}),
\begin{equation}\label{42}
\lim_{h\to\infty}\langle-{\rm div}(\text{\textbf{A}}_h\nabla u_h)\, ,\, v\rangle = \lim_{h\to\infty}\langle f_h\, ,\, v\rangle= \langle f\, ,\, v\rangle.
\end{equation}
By the uniqueness of weak limit, together with (\ref{41}) and (\ref{42}) we get
\begin{equation}\label{43}
-{\rm div}(\mathcal{M})=f\, .
\end{equation}
Since $\text{\textbf{A}}_h\xrightarrow{\;\,{\tiny G}\;\,}\text{\textbf{A}}$, by Theorem 3 it is also true that $\text{\textbf{A}}_h^t\xrightarrow{\;\,{\tiny G}\;\,}\text{\textbf{A}}^t$. Consider now
\begin{equation}\label{44}
\langle\text{\textbf{A}}_h\nabla u_h,\nabla v_h\rangle=\langle\nabla u_h,\text{\textbf{A}}_h^t\nabla v_h\rangle
\end{equation}
for a sequence $v_h\in H_0^1(\Omega)$ converging weakly to $v$ in $H_0^1(\Omega)$. The limit passage of (\ref{44}) together with Lemma 3 give
\begin{equation}\label{45}
\langle\mathcal{M},\nabla v\rangle=\langle\nabla u,A^t\nabla v\rangle\, ,
\end{equation}
hence
\begin{equation}\label{46}
\langle\mathcal{M},\nabla v\rangle=\langle\text{\textbf{A}}\nabla u,\nabla v\rangle\, .
\end{equation}
Take $\omega\subset\subset\Omega$ and $v_h$ such that $\nabla v=z\in \bb{R}^N$ on $\omega$, then (\ref{46}) can be written as
\begin{equation}\label{47}
\langle\mathcal{M}-\text{\textbf{A}}\nabla u,z\rangle=0,
\end{equation}
consequently, by the density of $v$ in $H_0^1(\Omega)$ we have
\begin{equation}\label{48}
\mathcal{M}-\text{\textbf{A}}\nabla u=0\, ,
\end{equation}
which completes the proof of the claim.

The following lemma proves the continuity of $F_h$ in $ H^1_0(\Omega)$.
\begin{Lem}
\emph{
$F_h(u)$ is continuous in $H^1_0(\Omega)$.
}
\end{Lem}
\hspace{-4mm}\underline{\emph{Proof}}. Let $u,v\in H^1_0(\Omega)$ be such that $\|u-v\|_{H^1_0(\Omega)}<\varepsilon$, then

\begin{eqnarray*}
\begin{array}{ll}
\vspace{1mm}
\hspace{-1.5cm}\big|F_h(u)-F_h(v)\big|\!\!\!&=\Big|\langle \mathbf{H}_hu,u\rangle - \langle \mathbf{H}_hv,v\rangle\Big|\\
&\leq\Big|\langle \mathbf{H}_h(u+v),(u-v)\rangle\Big|\\
&\leq\|\mathbf{H}_h(u+v)\|_{L^2(\Omega)}\,\|u-v\|_{L^2(\Omega)}\\
&\leq \varepsilon\|\mathbf{H}_h(u+v)\|_{L^2(\Omega)}.
\end{array}
\end{eqnarray*}
The term in the last inequality approaches zero as $\varepsilon\to0$, thus the lemma is proved.\hfill{$\blacksquare$}\\

The following theorem proves and characterizes the G-limit of $\mathbf{H}_h$ for another class of potentials $V_h$.
\begin{Theo}
\emph{
If $V_h$ is a weakly convergent sequence in $L^p(\Omega)$, $2\leq p<\infty$, with a weak limit denoted by $V$, then $\mathbf{H}_h$ G-converges to $\mathbf{H}=\mathbf{H}_0+V$.
}
\end{Theo}
\hspace{-4mm}\underline{\emph{Proof}}. Let $F_h$ and $F$ be the quadratic forms corresponding to $\mathbf{H}_h$ and $\mathbf{H}$ respectively. Following Theorem 4, to prove that $\mathbf{H}_h$ G-converges to $\mathbf{H}$, is equivalent to show that $F_h$ $\Gamma$-converges to $F$. To this end, we consider the quadratic form $F_h$ of $\mathbf{H}_h$,
\begin{equation*}
F_h(u)=\left\{ \begin{array}{ll}
\langle \mathbf{H}_hu,u\rangle\,,& u\in H^1_0(\Omega) \, , \\
\infty\,,& u\in L^2(\Omega)\backslash H^1_0(\Omega)\, .
\end{array} \right.
\end{equation*}
By the definition of $\Gamma$-convergence, to prove that $F$ is the $\Gamma$-limit of $F_h$, is equivalent to justify the following two conditions
\begin{itemize}
\item [$(i)$] $\liminf$-inequality: For every $u\in L^2(\Omega)$, and for every sequence $u_h$ converging to $u$ in $L^2(\Omega)$, $F(u)\leq\displaystyle\liminf_{h\to\infty} F_h(u_h)$.
\item [$(ii)$] $\lim$-equality: For every $u\in L^2(\Omega)$, there exists a sequence $u_h$ converging to $u$ in $L^2(\Omega)$ such that
    $F(u)=\displaystyle\lim_{h\to\infty} F_h(u_h)$.
\end{itemize}
To prove the $\liminf$-inequality we assume that $u_h\in H^1_0(\Omega)$. Otherwise the proof is obvious. By the continuity of $F_h$ in $H^1_0(\Omega)$ and since piecewise affine functions are dense in $H^1_0(\Omega)$, it suffices to prove the inequality for this class of functions (the same holds true for the $\lim$-equality).

Let $\Omega=\cup_{j=1}^m \Omega^{j}$ where $\Omega^{j}$ are disjoint sets, and let $u_h$ be linear in each $\Omega^{j}$ converging in $L^2(\Omega)$ to $u=\displaystyle\sum_{j=1}^m (a_jx+b_j)\chi_{\Omega^{j}}$, where $a_j$ and $b_j$ are elements of $\mathbb{C}^3$ and the product $a_jx$ is understood to be componentwise. Consider now $F_h$ with the sequence $u_h$,
\begin{equation}\label{83}
F_h(u_h)=\langle\mathbf{H}_0u_h,u_h\rangle + \langle V_hu_h,u_h\rangle\,.
\end{equation}
Since $u_h\to u$ in $L^2(\Omega)$,
\begin{equation}\label{84}
\langle\mathbf{H}_0u,u\rangle=\|\mathbf{H}^{1/2}_0u\|^2_{L^2(\Omega)}\leq\displaystyle \liminf_{h\to\infty} \|\mathbf{H}^{1/2}_0u_h\|^2_{L^2(\Omega)}=\displaystyle \liminf_{h\to\infty}\langle\mathbf{H}_0u_h,u_h \rangle\,.
\end{equation}
Hence
\begin{equation}\label{85}
\displaystyle\liminf_{h\to\infty}F_h(u_h)\geq \langle\mathbf{H}_0u,u\rangle\,+\,\displaystyle\liminf_{h\to\infty}\int_{\Omega}V_h|u_h|^2\,dx\,.
\end{equation}
For the last term of (\ref{85})
\begin{eqnarray}\label{2_2_2}
\begin{array}{lll}
\displaystyle\liminf_{h\to\infty}\displaystyle \int_{\Omega}V_h|u_h|^2\,dx & \!\!\!=&\!\!\!\!\!\!\!\!\!\displaystyle\liminf_{h\to\infty} \displaystyle\int_{\Omega} V_h|u+u_h-u|^2\,dx\\
&\!\!\!\geq&\!\!\!\!\!\!\!\displaystyle\liminf_{h\to\infty}\displaystyle \int_{\Omega}\!\!V_h|u|^2\,dx+\displaystyle \liminf_{h\to\infty}\displaystyle \int_{\Omega}\!\!V_h\,u^*\,(u_h-u)\,dx+\\
&\;+&\!\!\!\displaystyle\liminf_{h\to\infty}\displaystyle \int_{\Omega}\!\!V_h\,u\,(u_h-u)^*\,dx.
\end{array}
\end{eqnarray}
The symbol $*$ in (\ref{2_2_2}) refers to the complex conjugate. The first term to the right of the inequality of (\ref{2_2_2}) converges to $\displaystyle\int_{\Omega}V|u|^2\,dx$ by the weak convergence of $V_h$ to $V$ in $L^p(\Omega)$, $2\leq p<\infty$. Since $u_h\to u$ in $L^2(\Omega)$, the second and third terms to the right of the inequality of (\ref{2_2_2}) are vanishing as $h\to\infty$. Thus we have the $\liminf$-inequality, namely
\begin{equation}\label{86}
\displaystyle\liminf_{h\to\infty}F_h(u_h)\geq \langle\mathbf{H}_0u,u\rangle\,+\,\langle Vu,u\rangle=F(u).
\end{equation}

To prove the $\lim$-equality for some convergent sequence, again by the continuity argument it is enough to justify the equality for a piecewise affine sequence. So consider $u_h=u=(ax+b)\chi_{\Omega}$, then
\begin{equation*}
\begin{array}{ll}
\displaystyle\lim_{h\to\infty}F_h(u_h)&\!\!\!= \langle\mathbf{H}_0u,u\rangle\,+\,\displaystyle \lim_{h\to\infty}\langle V_hu,u\rangle\\
&\!\!\!=\langle\mathbf{H}_0u,u\rangle\,+\,\langle Vu,u\rangle\,,
\end{array}
\end{equation*}
the resulted limit is due to the boundedness of the set $\Omega$ and the linearity of $u$.\hfill{$\blacksquare$}\\

\end{document}